\documentclass{amsart}
\usepackage{amssymb,amsfonts,latexsym}

\newtheorem{theorem}{Theorem}[section]
\newtheorem{lemma}[theorem]{Lemma}
\newtheorem{proposition}[theorem]{Proposition}
\newtheorem{corollary}[theorem]{Corollary}
\theoremstyle{definition}

\newtheorem{remark}[theorem]{Remark}

\newcommand{\id}{\text{id}}

\newcommand{\End}{\text{End}}

\newcommand{\Hom}{\text{Hom}}
\newcommand{\Ad}{\text{Ad}}

\newcommand{\gr}{\text{gr}}

\newcommand{\Rep}{\text{Rep}}

\newcommand{\ot}{\otimes}

\newcommand{\ben}{\begin{enumerate}}
\newcommand{\een}{\end{enumerate}}

\newcommand{\Rad}{{\text{Rad}}}

\newcommand{\Q}{{\mathbb Q}}
\newcommand{\Z}{{\mathbb Z}}

\newcommand{\C}{{\mathbb C}}

\begin{document}

\title[basic Quasi-Hopf algebras of dimension $n^3$]
{basic Quasi-Hopf algebras of dimension $n^3$}
\author{Shlomo Gelaki}
\address{Department of Mathematics, Technion-Israel Institute of
Technology, Haifa 32000, Israel}
\email{gelaki@math.technion.ac.il}
\date{September 6, 2004}
\maketitle

\section{Introduction}

It is shown in \cite{eg} that a finite dimensional quasi-Hopf
algebra over $\mathbb{C}$ with radical of codimension $2$ is twist
equivalent to a Nichols Hopf algebra \cite{n}, or to a lifting of
one of four special quasi-Hopf algebras of dimensions 2, 8, 8, and
32.

The purpose of this paper is to construct new finite dimensional
basic \footnote{Recall that a finite dimensional algebra is called 
{\em basic} if every irreducible representation of it is $1-$dimensional.}
quasi-Hopf 
algebras $A(q)$ of dimension $n^3$, $n>2$,
parameterized by primitive roots of unity $q$ of order $n^2$, with
radical of codimension $n$, which generalize the construction of
the basic quasi-Hopf algebras of dimension $8$ given in \cite{eg}.
These quasi-Hopf algebras are not twist equivalent to a Hopf
algebra, and may be regarded as quasi-Hopf analogs of Taft Hopf
algebras \cite{t}.

By \cite{eo}, our construction is equivalent to the construction
of new finite tensor categories with integer Frobenius-Perron
dimensions of objects, whose simple objects form a cyclic group of
order $n$, and which are not tensor equivalent to a representation
category of a finite dimensional Hopf algebra.

We also prove that if $H$ is a finite dimensional radically graded
quasi-Hopf algebra with $H[0]=(k[\Z/n\Z],\Phi)$, where $n$ is
prime and $\Phi$ is a nontrivial associator, such that $H[1]$ is a
free left module over $H[0]$ of rank $1$ (it is always free by
\cite{s}) then $H$ is isomorphic to $A(q)$.

This paper is the first part of the author's work on nonsemisimple
quasi-Hopf algebras with a prime number of simple objects, all of
which are invertible. The second part of this work is the paper
math.QA/0403096 which is written jointly with Pavel Etingof.

{\bf Acknowledgments.} The author is grateful to P. Etingof for
numerous stimulating discussions.

The author's research was supported by The Israel
Science Foundation (grant No. 70/02-1).

\section{Preliminaries}

Let $H=(H,\Delta,\varepsilon,\Phi,S,\alpha,\beta)$ be a quasi-Hopf
algebra \cite{d}. Let $J=\sum_i f_i\ot g_i\in H\ot H$ be an
invertible element satisfying $$(\varepsilon\ot
\id)(J)=(\id\ot\varepsilon)(J)=1,$$ with inverse $J^{-1}=\sum_i
\overline{f_i}\ot \overline{g_i}$. Following \cite{d} we will call
such an element a {\em twist} for $H$. (In \cite{k}, it is called
a {\em gauge transformation}.) Set
$$\alpha_J:=\sum_iS(\overline{f_i})\alpha
\overline{g_i}\;\text{and}\;\beta_J:=\sum_i f_i\beta S(g_i).$$

In \cite{d} it is explained that given a twist $J$ for $H$, 
if $\beta_J$ is invertible then one
can define a new quasi-Hopf algebra structure
$H_J=(H,\Delta_J,\varepsilon,\Phi_J,S_J,\beta_J\alpha_J,1)$ on the
algebra $H$, where $$\Delta_J(h)=J\Delta(h)J^{-1},\;h\in H,$$
$$\Phi_J:=(1\ot J)(\id\ot \Delta)(J)\Phi(\Delta\ot
\id)(J^{-1})(J\ot 1)^{-1},$$ and $$S_J(h):=\beta_JS(h)\beta_J
^{-1},\;h\in H.$$

\begin{remark}
For basic quasi-Hopf algebras the elements
$\alpha$ and $\beta$ are always invertible. Indeed, if $\Phi=\sum_i X_i
\ot Y_i\ot Z_i$ then $\sum_i X_i\beta S(Y_i)\alpha Z_i=1$ (\cite{d}),
so $\alpha$ and $\beta$ cannot be zero on any $1-$dimensional representation 
and hence are invertible.
\end{remark}

\section{The main results}

\subsection{}
Let $n\ge 2$ be a positive integer, $k$ a field of characteristic
prime to $n$, and $q$ a fixed primitive root of unity of degree $n^2$.

Let $H=H(q)$ be the Taft Hopf algebra \cite{t} generated by $g,x$ with
relations
\begin{equation}\label{algrel}
x^{n^2}=0,\; g^{n^2}=1\; \text{and}\; gx=qxg,
\end{equation}
with coproduct
\begin{equation}\label{coalgrel}
\Delta(g)=g\otimes g,\; \Delta(x)=x\otimes g+1\otimes x,
\end{equation}
with counit
\begin{equation}\label{courel}
\varepsilon(g)=1,\; \varepsilon(x)=0,
\end{equation}
and with antipode
\begin{equation}\label{ant}
S(g)=g^{-1},\;S(x)=-xg^{-1}.
\end{equation}
Then $H$ is pointed (i.e. $H^*$ is basic) of dimension $n^4$.

Let $A=A(q)$ be the subalgebra of $H$ generated by $a:=g^n$ and $x$. Then
$\dim(A)=n^3$. Note that $A\subset H$ is not closed under
coproduct.

We will construct a twist $J$ in $k[\Z/n^2\Z]\ot k[\Z/n^2\Z]$ such
that the associator $\Phi:=dJ$ belongs to $k[\Z/n\Z]^{\ot
3}\subset k[\Z/n^2\Z]^{\ot 3}$ and $J\Delta(x)J^{-1}$ belongs to
$A\ot A$. By constructing such $J$ we will get that $A\subset H_J$
is a quasi-Hopf subalgebra.

\begin{remark}
The idea of constructing quasi-Hopf algebras as sub-quasi-Hopf algebras
of twists of Hopf algebras arose in a different context in \cite{er}.
\end{remark}

Let $\{1_z|0\le z\le n^2-1\}$ be the set of primitive idempotents
of $k[\Z/n^2\Z]$ defined by the condition $g1_z=q^z1_z$, and
denote by $i'$ the remainder of division of $i$ by $n$. Define
\begin{equation}\label{j}
J:=\sum_{z,y=0}^{n^2-1} c(z,y)1_z\otimes
1_y,\;\text{where}\;c(z,y):=q^{-z(y-y')}.
\end{equation}

\begin{lemma}\label{l1} $J$ is invertible and satisfies
$$(\varepsilon\ot \id)(J)=(\id\ot\varepsilon)(J)=1.$$
\end{lemma}

\begin{proof}
Clearly $J^{-1}=\sum_{z,y=0}^{n^2-1} c(z,y)^{-1}1_z\otimes 1_y$.

Since $(\varepsilon\ot \id)(J)=\sum_{y=0}^{n^2-1}c(0,y)1_y$ and
$(\id\ot\varepsilon)(J)=\sum_{z=0}^{n^2-1}c(z,0)1_z$, the second
claim follows.
\end{proof}

Let $\{\mathbf{1}_s|0\le s\le n-1\}$ be the set of primitive
idempotents of $k[\Z/n\Z]$. Note that
\begin{equation}\label{idem}
\sum_{i=0}^{n-1}1_{s+ni}=\mathbf{1}_{s}
\end{equation} for any
$0\le s\le n-1$.

Let $a$ be a generator of $\Z/n\Z$ and $f\in \Hom(\Z/n\Z, \Q/\Z)$
be given by $a\mapsto 1/n$. Lift $f$ to $\Q$ and set $F:=df\in
H^2(\Z/n\Z,\Z)$. The following lemma is well known.

\begin{lemma}
The map $H^1(\Z/n\Z,k^*)=Hom(\Z/n\Z,k^*)\to H^3(\Z/n\Z,k^*)$ given
by $\gamma\mapsto \gamma\cup F$ defines an isomorphism of groups.
\end{lemma}

Using this lemma we obtain the following.

\begin{lemma}
The map $\omega:(\Z/n\Z)^3\to k^*$ given by
$\omega(i,j,k)=q^{i(j+k-(j+k)')}$ is a nontrivial 3-cocycle, where
$i'$ denotes the remainder of the division of $i$ by $n$. In
particular, the element
\begin{equation}\label{assoc}
\Phi_l:=\sum_{i,j,k=0}^{p-1}q^{il(j+k-(j+k)')}{\bf1}_i\ot
{\bf1}_j\ot {\bf1}_k
\end{equation} is a nontrivial associator for $\Z/n\Z$ for any
$1\le l\le n-1$.

\end{lemma}

\begin{lemma}\label{l2}
The associator $\Phi_J$ equals $\Phi_{-1}$, so in particular
belongs to $A\ot A\ot A$ and is non-trivial.
\end{lemma}

\begin{proof}
Recall that $$\Phi_J=dJ=(1\ot J)(\id\ot \Delta)(J)(\Delta\ot
\id)(J^{-1})(J\ot 1)^{-1}.$$ We compute:
\begin{equation*}
\begin{split}
\Phi_J &=\sum_{z,y,w,u,s,t,p,q=0}^{n^2-1}
c(z,y)c(w,u)c(s,t)^{-1}c(p,q)^{-1} 1_w1_b1_p\ot
1_z1_a1_{s-b}1_q\ot 1_y1_{u-a}1_t\\ &= \sum_{z,y,w=0}^{n^2-1}
c(z,y)c(w,z+y)c(w+z,y)^{-1}c(w,z)^{-1} 1_w\ot 1_z\ot 1_y\\ &=
\sum_{z,y,w=0}^{n^2-1}q^{-w(z'+y'-(z+y)')}1_w\ot 1_z\ot 1_y\\
&=\sum_{z,y,w=0}^{n-1}q^{-w(z+y-(z+y)')}{\bf 1}_w\ot {\bf 1}_z\ot
{\bf 1}_y,
\end{split}
\end{equation*}
as desired.
\end{proof}

\begin{lemma}\label{l3}
The coproduct $\Delta_J(x)$ of $x$ in $H_J$ lies in $A\ot A$, and
is given by the following formula:
\begin{equation}\label{deltajx}
\Delta_J(x)=x\ot \sum_{y=0}^{n-1}q^y {\bf 1}_y + 1\ot (1-{\bf
1}_{0})x + a^{-1}\ot {\bf 1}_{0}x.
\end{equation}
\end{lemma}

\begin{proof}
We compute:
\begin{equation*}
\begin{split}
 \Delta_J(x)& = J\Delta(x)J^{-1} = J(x\ot g)J^{-1}+J(1\ot x)J^{-1}\\
 &= \sum_{z,y,a,b=0}^{n^2-1} c(z,y)c(a,b)^{-1}1_z x 1_a\ot 1_y g
 1_b+\sum_{z,y,a,b=0}^{n^2-1}c(z,y)c(a,b)^{-1}1_z 1_a\ot 1_y x
 1_b\\
 &= \sum_{z,y,a=0}^{n^2-1} c(z,y)c(a,y)^{-1}1_z x 1_a\ot 1_y g+
 \sum_{z,y,b=0}^{n^2-1}c(z,y)c(z,b)^{-1}1_z \ot 1_y x 1_b
\end{split}
\end{equation*}

Using the fact that $1_w x=x1_{w-1}$ for all $w$, and
$1_yg=q^y1_g$ for all $y$, we get that the last equation equals:
$$\sum_{z,y=0}^{n^2-1} q^y c(z,y)c(z-1,y)^{-1}x 1_{z-1}\ot 1_y +
\sum_{z,y=0}^{n^2-1}c(z,y)c(z,y-1)^{-1}1_z \ot x 1_{y-1},$$ which
equals $$x\ot \sum_{y=0}^{n-1}q^y {\bf 1}_y + 1\ot x(1-{\bf
1}_{n-1}) + a^{-1}\ot x{\bf 1}_{n-1},$$ or equivalently to $$x\ot
\sum_{y=0}^{n-1}q^y {\bf 1}_y + 1\ot (1-{\bf 1}_{0})x + a^{-1}\ot
{\bf 1}_{0}x,$$ as desired.

\end{proof}

\begin{proposition}\label{l4}
$(A,\Delta_J,\Phi_J)$ is a quasi-subbialgebra of $H_J$.
\end{proposition}

\begin{proof} The facts that $$(\id\ot \Delta_J)\Delta_J(x)=
\Phi_J(\Delta_J\ot\id)\Delta_J(x)\Phi_J^{-1},\;x\in A,$$ and
$$(\varepsilon\ot\id)\Delta_J(x)=(\id\ot\varepsilon)\Delta_J(x)=x,\;x\in
A,$$ follow since $H_J$ is a quasi-Hopf algebra. The result
follows now from Lemmas \ref{l2}, \ref{l3}.
\end{proof}

\begin{lemma}\label{betaj}
We have
\begin{equation}\label{betas}
\alpha_J\beta_J=a^{-1}\;\text{and}\;
S_J(x)=-x\sum_{z=0}^{n-1}q^{n-z}{\bf 1}_z.
\end{equation}
\end{lemma}

\begin{proof}
Let us first compute $\beta_J$. By definition,
$$\beta_J=\sum_{z,y=0}^{n^2-1}c(z,y)1_z1_{-y}=\sum_{z=0}
^{n^2-1}c(z,-z)1_z=\sum_{z=0}^{n^2-1}q^{n(\frac{z-z'}{n}+1)z}1_z=
\sum_{z=0}^{n^2-1}q^{(z-z'+n)z}1_{z}.$$

Similarly, $$\alpha_J=\sum_{z=0} ^{n^2-1}c(-z,z)^{-1}1_z
=\sum_{z=0}^{n^2-1}q^{-n(\frac{z-z'}{n})z}1_z=\sum_{z=0}^{n^2-1}q^{-(z-z')z}1_{z}.
$$
Therefore,
$$\alpha_J\beta_J=\sum_{z=0}^{n^2-1}q^{(z-z'+n)z-(z-z')z}1_{z}=\sum_{z=0}^{n^2-1}q^{nz}1_z=a^{-1},$$
as desired.

Clearly, $S_J(a)=a^{-1}\in A$. Let us now compute $S_J(x)$. By
definition, $S_J(x)=\beta_J(-xg^{-1})\beta_J^{-1}$. Since
$1_zx=x1_{z-1}$ and $1_zg^{-1}=q^{-z}1_{z}$, we get that
$$\beta_J(-xg^{-1})\beta_J^{-1}=-x\sum_{z=0}^{n^2-1}q^{(z+1-(z+1)'+n)(z+1)-z-(z-z'+n)z}1_z.$$
It is straightforward to check that the exponent of $q$ equals
$n-z'$ modulo $n^2$, and hence that
$$S_J(x)=-x\sum_{z=0}^{n^2-1}q^{n-z'}1_z=-x\sum_{z=0}^{n-1}q^{n-z}{\bf
1}_z,$$ as desired.
\end{proof}

Therefore, we have proved our first main result:

\begin{theorem}\label{main}
For any $q$,
$(A(q),\Delta_J,\varepsilon,\Phi_J,S_J,a^{-1},1)$ is a basic quasi-Hopf
algebra of dimension $n^3$, which is not twist equivalent to a Hopf algebra.
\end{theorem}

Theorem \ref{main} and \cite{eo} imply the following.

\begin{corollary}
The representation category $\Rep(A)$ of $A$ is a finite tensor
category with integer Frobenius-Perron dimensions of objects,
whose simple objects form a cyclic group of order $n$.
Furthermore, $\Rep(A)$ is {\em not} tensor equivalent to a
representation category of a finite dimensional Hopf algebra (i.e.
does not admit a fiber functor to the category of $k-$vector
spaces). However, it is a quotient of $\Rep(H)$ and has $(\Rep(\Z/n\Z),\Phi_J)$
as a semisimple
module category. So by \cite{eo}, it is tensor equivalent to a representation
category of a finite dimensional {\em weak} Hopf algebra.
\end{corollary}

\begin{remark}
In the case $n=2$, $A$ is isomorphic to $H_-(8)$ for $q=i$ and
$H_+(8)$ for $q=-i$, where $H_-(8)$, $H_+(8)$ are the
$8-$dimensional quasi-Hopf algebras constructed in \cite{eg}.
Indeed, one has to replace $x$ with $g^2x$.
\end{remark}

\subsection{}
Let $H$ be a finite dimensional radically graded quasi-Hopf algebra;
that is, $H=\gr H$ where $\gr H$ is taken with respect to the radical
filtration. Suppose that $H[0]=(k[\Z/n\Z],\Phi)$,
where $n>2$ is an integer and $\Phi$ is a nontrivial associator. (For
the case $n=2$, see \cite{eg}.)
In this case, as explained in \cite{eg}, $H$ is generated by $H[0]$ and $H[1]$.
By \cite{s}, $H[1]$ is a free left module over $H[0]$.
Assume that it is a module of {\bf rank $1$}.

We can assume that $\Phi=\Phi_s=
\sum_{i,j,k=0}^{n-1}Q^{si(j+k-(j+k)')/n}{\bf 1}_i\ot {\bf 1}_j\ot {\bf 1}_k,$
where $Q$ is a primitive root of $1$ of degree $n$.

Let $a$ be a generator of $\Z/n\Z$, and $Q^r$ be the scalar by which
of $\Ad(a)$ acts on $H[1]$.
By Theorem 2.17 in \cite{eo}, $r\ne 0$.

Assume that $r$ and $s$ are {\bf coprime} to $n$. Then by changing $Q$, we can
assume that $s=-1$,
and then, changing $Q$ to $Q^{r^2}$ and $a$ to $a^r$, we get to a situation
where $s=-1$ and $r=1$. Thus, we can assume from the beginning that $s=-1$ and
$r=1$.

Our second main result is the following.

\begin{theorem}\label{main2}
Let $H$ be as above. Then $H=A(q)$ for some primitive root of unity $q$
of degree $n^2$.
Moreover, the quasi-Hopf algebras $\{A(q)\}$ are pairwise nonisomorphic.
\end{theorem}

\begin{corollary}
Suppose $n$ is prime. Then $H=A(q)$ for some primitive root of unity $q$
of degree $n^2$.
\end{corollary}

\begin{proof}
For $n=2$ it is proved in \cite{eg}. For $n>2$ it follows from Theorem
\ref{main2}, and the fact that $r\ne 0$ by Theorem 2.17
in \cite{eo}.
\end{proof}

\begin{remark} Note that if the associator of $H$ is trivial (and $H[1]$ still
has rank $1$ over $H[0]=k[\Z/n\Z]$), then $H$ is obviously equivalent to the
Taft algebra of dimension $n^2$.
\end{remark}

The rest of the subsection is devoted to the proof of Theorem \ref{main2}.

Let $\chi$ be the character of $\Z/n\Z$ defined by $\chi(a)=Q$.
For $z\in H$ and $1\le l\le n$ define
\begin{equation}
\xi_l(z):=(\chi^l\otimes
\id)(\Delta(z))\;\text{and}\;\eta_l(z):=(\id\otimes
\chi^l)(\Delta(z)).
\end{equation}
Then $\xi_l$ and $\eta_l$ are algebra homomorphisms
from $H$ to $H$, $\xi_l(a)=\eta_l(a)=Q^la$ and
$H[1]$ is invariant under $\xi_l$ and $\eta_l$.

Let $L_a:H\to H$ denote left multiplication by $a$. We will abuse notation
and write $a$ for $L_a$.

\begin{lemma}
$\xi_la=Q^la\xi_l$ and $\eta_la=Q^la\eta_l$.
\end{lemma}

Set
\begin{equation}
E_r=E_r(Q):=\sum_{k=0}^{n-2}{\bf 1}_{k-r}+Q{\bf 1}_{n-1-r},
\end{equation}
and $\xi:=\xi_1$, $\eta:=\eta_1$.

\begin{lemma}
$\xi^n=Q^{-1}\id$, $\eta^n=Q\id$ and $\xi\eta=E_0^{-1}E_{-1}\eta\xi$.
\end{lemma}

\begin{proof}
We have $$\xi(\xi_l(z))=(\chi^l\ot \chi\ot \id)((\id\ot \Delta)\Delta(z))=
(\chi^l\ot \chi\ot \id)(\Phi(\Delta\ot \id)\Delta(z)\Phi^{-1}).$$ Since
$(\chi^l\ot \chi\otimes \id)(\Phi)=\sum_{k=0}^{n-1}Q^{-l\frac{1+k-(1+k)'}{n}}{\bf
1}_k=\sum_{k=0}^{n-2}{\bf 1}_k+Q^{-l}{\bf 1}_{n-1}=E_0^{-l}$, we get that
$\xi(\xi_l(z))=E_0^{-l}\xi_{l+1}E_0^{l}=E_0^{-l}E_{-1}^{l}\xi_{l+1}$ and hence,
$\xi_{l+1}=E_0^{l}E_{-1}^{-l}\xi\xi_l$. So by induction, $$\xi_{m+1}=
E_0^{-m}E_{-1}E_{-2}\cdots E_{-m}\xi^{m+1}.$$ Therefore,
$\id=\xi_n=E_0E _{-1}E_{-2}\cdots E_{-(n-1)}\xi^{n}$. Hence,
we get that $\xi^n=Q^{-1}\id$, as desired.

Similarly, since $(\id\otimes \chi\ot \chi^l)(\Phi)$ equals $1$ if $l<n-1$
and $a$ if $l=n-1$, we get that $\eta^n=Q\id$.

Finally, $$\xi\eta(z)=(\chi\ot
\id \ot \chi)(\Phi(\Delta\ot \id)\Delta(z)\Phi^{-1})=
E_0^{-1}(\chi\ot \id \ot \chi)((\Delta\ot \id)\Delta(z))E_0,$$
which is equivalent to $\xi\eta=E_0^{-1}E_{-1}\eta\xi$, as desired.
\end{proof}

Let $B(Q)$ be the algebra generated by
$\xi,\eta,a$ with the following relations:
\begin{equation}
a^n=1,\; \xi^n=Q^{-1},\; \eta^n=Q,\;
\xi a=Qa\xi,\; \eta a=Qa\eta\;\text{and}\;\xi\eta=E_0^{-1}E_{-1}\eta\xi.
\end{equation}
We can rescale the generators so that $B(Q)$ is generated by
$\xi,\eta,a$ with the same relations except that
\begin{equation}
\xi^n=\eta^n=1.
\end{equation}
We have thus proved the following proposition.

\begin{proposition}\label{acts}
The algebra $B(Q)$ acts on the space $H[1]$.
\end{proposition}

Let us now study the representations of the algebra $B(Q)$.

Consider the $n$ quasi-Hopf algebras $A(q)$, where $q$ is
a $n-$th root of $Q$. In it $\Delta(x)$ is given by equation
(\ref{deltajx}). Set $V_{q}$
to be $A(q)[1]$. Then $V_{q}$ has a basis $\{{\bf 1}_ix|0\le i\le n-1\}$,
and the algebra $B(Q)$ acts on $V_{q}$ via $$a{\bf 1}_ix=Q^i{\bf 1}_ix\;,
\eta({\bf 1}_ix)=q{\bf 1}_{i-1}x\;\text{and}\;\xi({\bf 1}_ix)=
Q^{-\delta_{i,0}}{\bf 1}_{i-1}x.$$
In particular, the eigenvalues of $\eta\xi^{-1}$ on
$V_{q}$ are $q$ with
multiplicity $n-1$ and $Qq$ with multiplicity $1$.

\begin{proposition}\label{ss}
$B(Q)\cong \bigoplus_{q,\;q ^n=Q}\End_{\C}(V_{q})$.
\end{proposition}

\begin{proof}
Since the eigenvalues of $\eta\xi^{-1}$ are different for different choices
of $q$, the representations $V_{q}$ are pairwise non-isomorphic.
Also, the representations $V_{q}$ are irreducible. Therefore, there is a
surjective map $\varphi:B(Q)\to \bigoplus_{q,\;q ^n=Q}
\End_{\C}(V_{q})$. But, $\{a^i\xi^j\eta^k|0\le i,j,k\le n-1\}$ is a
spanning set for $B(Q)$, so the dimension of $B(Q)$ is less than or equal
to $n^3$. Therefore, $\varphi$
is an isomorphism, and so $B(Q)$ is semisimple, has dimension $n^3$ and
its irreducible representations are the $V_{q}$.
\end{proof}

\begin{corollary}
The quasi-Hopf algebras $A(q)$ are pairwise non-isomorphic.
\end{corollary}

\begin{proof}
Let $f: A(q_1)\to A(q_2)_K$ be an isomorphism of quasi-Hopf
algebras, where $K$ is a twist.
By passing to $\gr(f)$ (under the radical filtration), we may assume that
$K$ has degree $0$ and $f$ preserves the gradings. Then $f(a)=a^r$, where
$r$ is coprime to $n$.
In this case the condition that
$f$ maps the associator of $A(q_1)$
to the associator of $A(q_2)$ twisted by $K$ is that $q_1^n=q_2^{nr^2}$,
and the fact that the action of $Ad(a)$ on $A(q_i)[1]$, $i=1,2$, is by
$q_i^n$ implies
$q_1^n=q_2^{nr}$. This implies that $r=1$ (i.e. $q_1^n=q_2^n=Q$)
and $f(a)=a$. In particular, $K$ is a Hopf twist in $k[\Z/n\Z]^{\otimes
2}$, hence $K$ is a coboundary and hence we may assume that $K=1\otimes 1$.
So $f$ defines an isomorphism of representations of $B(Q)$:
$V_{q_1}\to V_{q_2}$. Hence $q_1=q_2$.
\end{proof}

Finally, by Propositions \ref{acts}, \ref{ss}, and the assumption that
$H[1]$ has rank $1$,
we have that $H[1]=V_q$ for some $q$. Let $v_1,\dots,v_n$ be the basis
of $V_q$ as above. Let $x:=\sum_j v_j$.
Then the coproduct formula for $x$ is as in (\ref{deltajx})
(since $\xi_j(x)$ and $\eta_j(x)$ are as in $A(q)$). Thus,
$H$ is a quotient of $A(q)$ by a graded ideal for the relevant $q$.

\begin{lemma}
$A=A(q)$ has no nonzero quasi-Hopf ideals other than itself and $I:=\Rad(A)$.
\end{lemma}

\begin{proof}
Let $L$ be a nonzero quasi-Hopf ideal of $A$, not equal to $A$ or $\Rad(A)$. Consider
$A':=\gr(A/L)$. Then $A'$ is a quasi-Hopf algebra and $A'[1]$ is a free
rank $1$ module over $A'[0]=(k[G],\Phi)$. Let $H'$ be obtained from $A'$ by
adjoining a new grouplike element $h$ such that $\Ad(h)$ acts by $q$ on
$A'[1]$ and $h^n=a$. Let $H'':=(H')_{J^{-1}}$, where $J$ is as in (\ref{j}).
Then $H''$ is a Hopf algebra generated by a grouplike element $h$
and a nonzero skew-primitive element $x$ satisfying
the relations of the Taft algebra of dimension $n^4$.
So $H''$ is a nonsemisimple quotient of the Taft algebra.
But it is well known that the Taft algebra does not have nontrivial
nonsemisimple quotients. So $H''$ is the Taft algebra, hence $\dim(A')=n^3$
and $A'=A$. Thus, $L=0$.
\end{proof}

We conclude that $H=A(q)$. This completes the proof of Theorem \ref{main2}.

\end{document}